\documentclass[12pt,a4paper]{article}
\usepackage{amsmath,amssymb,amsfonts}
\def\Bbb{\mathbb}

\title{\bf Chebyshev polynomials and Galois groups of De Moivre polynomials}

\author{Kurt Girstmair}

\date{}

\makeatletter
\let\@@maketitle=\maketitle
\def\maketitle{\def\thispagestyle##1{\relax}\@@maketitle}
\makeatother
\newtheorem{theorem}{Theorem}
\newtheorem{prop}{Proposition}
\newtheorem{lemma}{Lemma}
\newtheorem{corollary}{Corollary}

\def\BE{\begin{equation}}
\def\EE{\end{equation}}
\def\BD{\begin{displaymath}}
\def\ED{\end{displaymath}}
\def\BA{\begin{array}}
\def\EA{\end{array}}
\def\BEA{\begin{eqnarray*}}
\def\EEA{\end{eqnarray*}}
\def\BI{\bibitem}

\def\Z{\Bbb Z}
\def\Q{\Bbb Q}

\def\C{\Bbb C}

\def\OO{{\cal O}}

\def\phi{\varphi}

\def\MB{\mbox}
\def\LD{\ldots}

\def\SP#1{\langle #1 \rangle}

\def\sminus{\smallsetminus}
\def\DIV{\,|\,}
\def\NDIV{\, \nmid \,}

\def\BQ{``}
\def\EQ{'' }

\def\NI{\noindent}
\def\MN{\medskip\noindent}
\def\STOP{\hfill$\Box$}

\def\GAL{{\rm{Gal}}}

\begin{document}
\maketitle

\begin{abstract}
\NI
Let $n\ge 3$ be an odd natural number. In 1738, Abraham de Moivre introduced a family of polynomials of degree $n$ with rational coefficients, all of which are solvable.
So far, the Galois groups of these polynomials
have been investigated
only for prime numbers $n$ and under special assumptions. We describe the Galois groups for arbitrary odd $n\ge 3$ in the irreducible case, up to few exceptions.
In addition, we express all zeros of such a polynomial as rational functions of three zeros, two of which are connected in a certain sense. These results are based on the reduction of the radical
\BD
  \sqrt[n]{d+\sqrt R},
\ED
whose degree is $2n$ in general, to irrationals of degree $\le n$. Such a reduction was given in a previous paper of the author.
Here, however, we present a much simpler approach that is based on properties of Chebyshev polynomials.
\end{abstract}

\section*{1. Introduction}

Chebyshev polynomials were introduced by Abraham de Moivre in 1738, more than eighty years before Chebyshev was born. Indeed, De Moivre defined the polynomial $T_n\in\Q[Z]$ of degree $n$ by the
identity
\BE
\label{1.2}
  \cos(nx)=T_n(\cos(x))
\EE
and noted the formula
\BE
\label{1.4}
  T_n=\sum_{k=0}^{\lfloor n/2\rfloor}(-1)^k\frac{n}{n-k}\binom{n-k}{k}2^{n-2k-1}Z^{n-2k}
\EE
for $n\ge 1$, see \cite[p. 246]{Sch}.  Here $\binom{n-k}{k}$ is the usual binomial coefficient. From (\ref{1.2}) we see that $T_0=1$. Both (\ref{1.2}) and (\ref{1.4})
can be found in standard texts about Chebyshev polynomials, see \cite[p. 1, p. 32]{Ri}.

De Moivre used Chebyshev polynomials in the following context. Let $d$, $v$ be rational numbers, $d\ne 0$, and $R=d^2-v^n$ not a square.
We denote one of the square-roots of $R$ by $\sqrt R$ and the other by $-\sqrt R$.
Let
\BE
\label{1.9}
  h=Z^n-(d+\sqrt R),\: h'=Z^n-(d-\sqrt R).
\EE
Moreover, let $y,y'\in\C$ be zeros of $h,h'$, respectively. Since $(yy')^n=d^2-R=v^n$, these zeros can be chosen such that $yy'=v$.
Put $2x=y+y'$. Then $x$ satisfies the equation
\BD
\label{1.8}
  \sqrt v^{\:n}T_n(x/\sqrt v)-d=0,
\ED
see \cite[p. 246]{Sch}.
Therefore, we call $\sqrt v^{\: n}T_n(Z/\sqrt v)-d$ the {\em De Moivre polynomial} of degree $n$ with parameters $d,v$. It is clear that this polynomial has rational coefficients.
By (\ref{1.9}), the De Moivre polynomial is a solvable polynomial of degree $n$.

In the sequel we avoid the powers of 2 in formula (\ref{1.4})
and work with
\BE
\label{1.14}
  F_n=2T_n(Z/2)
\EE
instead of $T_n$. The corresponding De Moivre polynomial takes the form
\BE
\label{1.16}
  \sqrt v^{\:n}F_n(Z/\sqrt v)-2d
\EE
and has the zero $2x$. The polynomial $\sqrt v^{\:n}F_n(Z/\sqrt v)$ is also known as the Dickson polynomial of degree $n$.

If $n$ is odd, the following kind of De Moivre polynomial is more appropriate for our purpose.
Let $d,R$ be rational numbers, $d\ne 0$, $R$ not a square. Form the polynomials $h$, $h'$ as above, see (\ref{1.9}). Let $y,y'$ be zeros of $h,h'$, respectively.
Put
\BE
\label{1.10}
   z=yy'\: \MB{ and }\: u=z^{(n-1)/2}(y+y').
\EE
Then $u$ is a zero of
\BE
\label{1.12}
 f_n= \sqrt D^{\:n}F_n(Z/\sqrt D)-dD^{(n-1)/2},
\EE
where $D=d^2-R$. In contrast with the general case, $R$ need not have a special form and $y,y'$ are not subject to an additional condition.
Obviously, $f_n$ is the De Moivre polynomial of degree $n$ (in the sense of (\ref{1.16})) with the parameters $v'=D$ and $d'=dD^{(n-1)/2}$.

Suppose that $n\ge 3$ is odd.
The main purpose of this paper is the study of De Moivre polynomials $f_n$ as in (\ref{1.12}), for which we also write $f_n(Z,d,R)$ (recall $D=d^2-R$).

In the literature the Galois group of the splitting field $L$ of $f_n$ has been studied only in certain cases for a prime numbers $n=p\ge 3$.
Indeed, for $p=5$, $\GAL(L/\Q)$ has been determined by Borger \cite{Bo}, and Spearman and Williams \cite{SpWi}.
Bruen, Jensen, and Yui consider the case $p$ prime, $p\equiv 3$ mod 4 and $R=s^2(-p)$, $s\in\Q$, see \cite{Br}. They find that $|\GAL(L/\Q)|=p(p-1)/2$,
provided that $f_n=f_p$ is irreducible. Filaseta, Luca, St\u anica, and Underwood study the case $d=-1/(2m^{(p-1)/2})$ and $R=(1-4m^p)/(4m^{p-1})$ (where $D$ equals $m$) for natural numbers $m\ge 2$.
They show that $f_p$ is irreducible and $|\GAL(L/\Q)|=p(p-1)$, see \cite{Fi}.

In what follows we write $R=s^2R'$, where $s$ is a rational number and $R'\in\Z$ is square-free. Moreover, let $\zeta_n$ denote a primitive $n$th root of unity in $\C$.
The main result of the present paper is the following.

\begin{theorem} 
\label{t1}
Let $n\ge 3$ be an odd natural number, $d,R\in\Q$, $d\ne 0$, $R$ not a square, and $D=d^2-R$. Let $L\subseteq\C$ be the splitting field of $f_n=f_n(Z,d,R)$.
Suppose, further, that $f_n$ is irreducible in $\Q[Z]$ and that $R'\ne -3$ if $\:3\DIV n$. Then
$\GAL(L/\Q)$ is isomorphic to a semidirect product
\BD
  (\Z/n\Z)\rtimes (\Z/n\Z)^{\times}
\ED
except in the case $R'<0$ and $\sqrt{R'}\in\Q(\zeta_n)$, where $\GAL(L/\Q)$ is isomorphic to
\BD
 (\Z/n\Z)\rtimes ((\Z/n\Z)^{\times}/\{\pm 1\}).
\ED
\end{theorem} 

\MN
Theorem \ref{t1} also holds in the case $3\DIV n$ and $R'=-3$ under further restrictions, see Section 4.
In the paper \cite{Br} a necessary and sufficient condition for the irreducibility of $f_n$ is given for the special case considered there. The following result both simplifies and generalizes
this criterion.

\begin{theorem} 
\label{t2}
In the setting of Theorem \ref{t1}, the polynomial $f_n$ is irreducible if, and only if, for all primes $p$ dividing $n$, $(d+\sqrt R)/(d-\sqrt R)$ is not a $p$th power in $\Q(\sqrt R)$.
\end{theorem} 

\MN
Theorem \ref{t2} implies \BQ rational\EQ criteria for the irreducibility of $f_n$ that are easy to check, see Proposition \ref{p7}.

One of our basic tools is an expression for the above zeros $y$, $y'$ of $h,h'$, respectively, as a polynomial function of $u$, $z$ and $\sqrt R$ (recall (\ref{1.9}) and (\ref{1.10})). In the previous paper \cite{Gi} we exhibited
the respective polynomial.
The proof was based on hypergeometric summation, in particular, on Zeilberger's algorithm. In the meanwhile, however, we learned that all polynomials involved can be expressed in terms of Chebyshev polynomials.
Thereby, one can prove the result of \cite{Gi} by means of well known identities for Chebyshev polynomials. This will be the subject of Section 2.

Section 3 describes the splitting field of $f_n$ over $\Q$. Moreover, it contains a formula that expresses all zeros of $f_n$ in terms of three of them, two of which, however, are connected in a certain sense
(Theorem \ref{t3}).

In Section 4 we prove the above Theorems. We also show that if $n=p$ is a prime, then either $f_p$ is irreducible or has a rational zero (Theorem \ref{t4}). This result allows proving the irrationality of
the polynomial $f_p$ in the case considered by Filaseta et al. \cite{Fi} in a simple way.

Section 5 consists of some remarks on the case of an even degree $n$.

Finally, we mention two more results. Abhyankar, Cohen and Zieve determine the Galois group of the splitting field of $\sqrt v^{\:n}F_n(Z/\sqrt v)-X$ over the rational function field $K(X)$ under the assumption that $v\ne 0$
and that the ground field
$K$ of a characteristic not dividing $n$ contains a primitive $n$th root of unity. They show that this group is a dihedral group of order $2n$ if $n\ge 3$, see \cite{Abh}.
In the paper \cite{Ki} Kida considers the case that $n$ is odd and square-free. He shows that every cyclic extension of degree $n$ of the maximal real subfield $\Q(\zeta_n)^+$
of $\Q(\zeta_n)$ is generated by a polynomial of the form $T_n-a$, where $a\in\Q(\zeta_n)^+$ is defined explicitly.

\section*{2. Reduction by means of Chebyshev polynomials}

Let $n\ge 3$ be an odd natural number, $d,R\in\Q$, $d\ne 0$, $R$ not a square. As above, let $y$, $y'$ be zeros of the polynomials $h$, $h'$, respectively (see (\ref{1.9})). In general, $y$ and $y'$ are
irrationals of degree $2n$ over $\Q$. In this section we reduce these irrationals to irrationals of degree $\le n$ over $\Q$, i.e., we express the former as a polynomial function of the latter. To this end we
define the polynomial
\BE
\label{2.0}
 A=\left(F_{n-1}(Z/\sqrt D)-\frac{dZ}D\right)\frac 1{2R}\in\Q[Z]
\EE
(recall (\ref{1.14}) and the fact that $n-1$ is even). Further, we put $z=yy'$ and $u=z^{(n-1)/2}(y+y')$, see (\ref{1.10}).
Then $z$ is a zero of $Z^n-D\in\Q[Z]$ and $u$ a zero of $f_n=f_n(Z,d,R)$. In the paper \cite{Gi} we proved the following result.

\begin{prop} 
\label{p1}
In the above setting,
\BE
\label{2.2}
  \{y,y'\}=\left\{z^{(n+1)/2}\left( \frac u{2D} \pm A(u)\sqrt{R}\right)\right\}.
\EE
\end{prop} 

\MN
As mentioned above, our proof of (\ref{2.2}) was based on hypergeometric summation. Here we give a simple proof in terms of well-known identities of Chebyshev polynomials.
To this end we put
\BE
\label{2.3}
 f_{n-2}'=\left(\frac 1{\sqrt D^{\:n-4}}F_{n-2}(Z/\sqrt D)-\frac{2d}{D^{(n-3)/2}}\right)\frac 1R.
\EE
We are going to prove the relation
\BE
\label{2.4}
4D^2RA^2=f_nf_{n-2}'+Z^2-4D.
\EE
Suppose, for the moment, that (\ref{2.4}) holds. Since $u$ is a zero of $f_n$, we have $4D^2RA^2(u)=u^2-4D=z^{n-1}((y+y')^2-4yy')$. This gives $2DA(u)\sqrt R=\pm z^{(n-1)/2}(y-y')$.
If we observe $u=z^{(n-1)/2}(y+y')$, we readily obtain (\ref{2.2}).

In order to prove (\ref{2.4}), we note the relation
\BD
 T_nT_m=\frac 12(T_{n+m}+T_{n-m})
\ED
for Chebyshev polynomials with $n\ge m$, see \cite[p. 5]{Ri}. Since $T_0=1$, $T_1=Z$, and $T_2=2Z^2-1$, this relation implies
\BD
  T_{n-1}^2=T_nT_{n-2}-Z^2+1.
\ED
From (\ref{1.14}) we obtain
\BE
\label{2.6}
  F_{n-1}^2=F_nF_{n-2}-Z^2+4.
\EE
We also note
\BD
 T_n=2ZT_{n-1}-T_{n-2},
\ED
see \cite[p. 35]{Ri}. This is equivalent to
\BE
\label{2.8}
 ZF_{n-1}=F_n+F_{n-2}.
\EE
Now we expand both sides of (\ref{2.4}) by means of (\ref{2.0}), (\ref{1.12}), and (\ref{2.3}). Furthermore,  we use the identities (\ref{2.6}) and (\ref{2.8}). In this way we obtain the desired result.

\section*{3. The splitting field and a formula for the zeros}

Throughout this section, let $n\ge 3$ be odd, $d,R\in\Q$, $d\ne 0$,  $R$ not a square, $D=d^2-R$, and $f_n=f_n(Z,d,R)$ (see (\ref{1.12})). Let $y,y'$ be zeros of $h$, $h'$, respectively (see (\ref{1.9})), $z=yy'$
and $u=z^{(n-1)/2}(y+y')$ (see (\ref{1.10})). Then $u$ is a zero of $f_n=f_n(Z,d,R)$. Let $\zeta=\zeta_n$ denote a primitive $n$th root of unity.
From the context of formulas (\ref{1.9}) -- (\ref{1.12}) we conclude that the numbers
\BE
\label{3.2}
  u_k=z^{(n-1)/2}(y\zeta^k+y'\zeta^{-k}),
\EE
$k=0,\LD,n-1$, are zeros of $f_n$ (with $u=u_0$) since $y\zeta^k$ is a zero of $h$ and $y'\zeta^{-k}$ a zero of $h'$.  These numbers are all distinct. Indeed,
if $u_k=u_l$, then $y\zeta^k$ and $y\zeta^l$ are both solutions of the quadratic equation
\BD
  z^{(n-1)/2}(x+z/x)=u_k
\ED
in $x$. However, this equation has exactly two solutions, namely, $y\zeta^k$ and $y'\zeta^{-k}$. Since $y\ne 0$, we have $y\zeta^l=y'\zeta^{-k}$ and $(y\zeta^l)^n=d+\sqrt R=(y'\zeta^{-k})^n=d-\sqrt R$.
Therefore, $R=0$, a contradiction. Accordingly, $f_n$ has the distinct zeros $u_k$, $k=0,\LD, n-1$, and $L=\Q(u, u_1, \LD, u_{n-1})$ is the splitting field of $f_n$ over $\Q$.

\MN
{\em Remark.}
If we choose a zero of $f$ different from $u$, say $u'=u_j$, and form the zeros of $f$ according
to (\ref{3.2}) (with $y\zeta^j$ instead of $y$ and $y'\zeta^{-j}$ instead of $y'$), we obtain
$u_k'=u_{k+j}$. Here the indices $j$, $k$, and $k+j$ must be interpreted as elements of $\Z/n\Z$. We will use this convention a number of times in this paper.

\MN
For the sake of simplicity we assume, henceforth, that the plus sign in (\ref{2.2}) gives $y$ and the minus sign gives $y'$.
In the opposite case, we interchange $\zeta$ and $\zeta^{-1}$ to the effect that the plus sign in (\ref{2.2}) yields $y$.
Then (\ref{2.4}) implies the following formula:
\BE
\label{3.4}
   u_k=\frac u2(\zeta^k+\zeta^{-k})+ DA(u)\sqrt R(\zeta^k-\zeta^{-k}),
\EE
$k=0,\LD,n-1$.

\MN
{\em Remark.}
Since $n$ is odd, the polynomial $f_n\in\Q[Z]$ has a real zero $u$. Then (\ref{3.4}) shows that all other zeros of $f_n$ are pairwise complex-conjugate if $R>0$. It also shows that all
zeros of $f$ are real if $R<0$.

\MN
Let $K=\Q(\sqrt R(\zeta-\zeta^{-1}))$. We have the following.

\begin{prop} 
\label{p2}
In the above setting, let $u$ be a zero of $f_n$. Then the splitting field $L$ of $f_n$ is
the composite field
\BD
 \Q(u)\cdot K.
\ED
\end{prop} 

\MN
{\em Proof.}
In view of the above remark, we may assume that $u=u_0=z^{(n-1)/2}(y+y')$ and that the zeros $u_k$ of $f$ are given by (\ref{3.4}).
We observe $(u_1-u_{n-1})/2=DA(u)\sqrt R(\zeta-\zeta^{-1})\in L$. But $D\ne 0$, and since $y\ne y'$, formula (\ref{2.2}) shows $A(u)\ne 0$.
Hence $\sqrt R(\zeta-\zeta^{-1})\in L$ and $\Q(u)\cdot K\subseteq L$.

Conversely, we observe that $K$ is a subfield of the abelian extension $\Q(\sqrt R,\zeta)$ of $\Q$. Since
$(\sqrt R(\zeta-\zeta^{-1}))^2=R(\zeta^2+\zeta^{-2}-2)$, we have $\zeta^2+\zeta^{-2}\in K$. But $\zeta^2+\zeta^{-2}$ generates the maximal real subfield
$\Q(\zeta)^+=\Q(\zeta+\zeta^{-1})$ of the $n$th cyclotomic field $\Q(\zeta)$. Since $K$ is a Galois extension of $\Q$, we have $\Q(\zeta)^+\subseteq K$.
The numbers $\zeta^k+\zeta^{-k}$, $k\in\{0,\LD,n-1\}$, however, are real,
and so they also lie in $\Q(\zeta)^+$ and, thus, in $K$.

Next let $k\in\{1,\LD,n-1\}$ be odd. Then (\ref{1.2}) and (\ref{1.14}) give
\BE
\label{3.6}
  2\cos(kx)=F_k(2\cos(x)).
\EE
If we apply this formula to $\pi/2-x$ instead of $x$, we obtain
\BE
\label{3.7}
 2\sin(kx)=(-1)^{(k-1)/2}F_k(2\sin(x)).
\EE
Note that the polynomial
\BE
\label{3.8}
 P_k=i\sqrt R (-1)^{(k-1)/2}F_k(Z/(i\sqrt R))
\EE
has rational coefficients.
Since $\sqrt R(\zeta^k-\zeta^{-k})=i\sqrt R\cdot 2\sin(2\pi k/n)$ and $\sqrt R(\zeta-\zeta^{-1})=i\sqrt R\cdot 2\sin(2\pi/n)$,
the identity (\ref{3.7}) gives
\BE
\label{3.10}
\sqrt R(\zeta^k-\zeta^{-k})=P_k(\sqrt R(\zeta-\zeta^{-1})).
\EE
We also observe $\sqrt R(\zeta^{n-k}-\zeta^{-(n-k)})=-\sqrt R(\zeta^k-\zeta^{-k})$. Accordingly, $\sqrt R(\zeta^k-\zeta^{-k})\in K$ for all $k\in\{0,\LD,n-1\}$.
By (\ref{3.4}), all zeros $u$ of $f_n$ are in $\Q(u)\cdot K$.
\STOP

\begin{corollary} 
\label{c1}
If $f_n$ has a zero $u\in\Q$, then $L=K$.

\end{corollary} 

\MN
Now we are in a position to express all zeros of $f_n=f_n(Z,d,R)$ in terms of three zeros, two of which, however, are connected. Recall that $n$ is odd, $d,R\in\Q$, $d\ne 0$. and $R$ not a square.

\begin{theorem} 
\label{t3}
Let $u$ be a zero of $f_n=f_n(Z,d,R)$ and $\zeta$ a primitive $n$th root of unity. Define the zeros $u_1$ and $u_{n-1}$ by formula {\rm (\ref{3.4})} (with $k=1$ and $k=n-1$).
Let $k\in\{1,\LD,n-1\}$ be odd and $F_k$ and $P_k$ as in {\rm (\ref{1.14})} and {\rm (\ref{3.8})}. Then
\BD
\label{3.12}
  u_k=\frac u2 F_k\left(\frac{u_1+u_{n-1}}{u}\right)+ DA(u)P_k\left(\frac{u_1-u_{n-1}}{2DA(u)}\right).
\ED
We obtain $u_{n-k}$ if we replace the plus sign in the middle of this formula by the minus sign.
\end{theorem} 

\MN
{\em Proof.} Formula (\ref{3.6}) yields
\BE
\label{3.14}
  \zeta^k+\zeta^{-k}=F_k(\zeta+\zeta^{-1}).
\EE
By (\ref{3.4}), $\zeta+\zeta^{-1}=(u_1+u_{n-1})/u$ and $\sqrt R(\zeta-\zeta^{-1})=(u_1-u_{n-1})/(2DA(u))$ (recall $A(u)\ne 0$). Therefore, (\ref{3.4}), (\ref{3.14}), and (\ref{3.10}) prove the result.
\STOP

\MN
{\em Example.} Let $n=9$, $d=26$, $R=675=15^2 3$. The polynomial $f_9$ is reducible and has the decomposition $f_9=(Z^3-3Z-4)(Z^6-6Z^4+4Z^3+9Z^2-12Z+13)$. Let $u$ be the real zero of the first factor and
$\zeta=e^{2\pi i/9}$. Then $u_1\approx 1.682098-0.582651i$ is a zero of the second factor and $u_8$ is the complex-conjugate of $u_1$. Moreover, $A(u)\approx -0.017445$. 
For $k=7$, we have $F_k=Z^7-7Z^5+14Z^3-7Z$ and $P_k=(1/R^3)Z^7+(7/R^2)Z^5+(14/R)Z^3+7Z$.
Theorem \ref{t3} yields  $u_7\approx 0.381301+0.892673i$, another zero of the second factor.

\MN
{\em Remark.}
We briefly highlight Theorem \ref{t3} in the case when $n=p$ is a prime. Here all $\zeta^k$, $k\in\{1,\LD,p-1\}$, are primitive $p$th roots of unity. Therefore, $u_1$ may be an arbitrary zero $\ne u$ of $f_p$, i.e.,
$\zeta$ may be defined by (\ref{3.4}). But then $u_{p-1}$ must satisfy (\ref{3.4}) with $k=p-1$.
A theorem of Galois says that an irreducible polynomial of degree $p$ is solvable if, and only if, all of its zeros can be written as rational functions of two of them, where the
latter may be chosen arbitrarily, see \cite[p. 163]{Hu}. In the paper \cite{SpWi} such
functions have been given for $f_p$ in the case $p=5$. Here we have a slightly weaker result for all $p$. If $R>0$ and $u$ is the real zero of $f$, ($u_1,u_{p-1})$ may be chosen
as an arbitrary pair of complex-conjugate zeros (recall the remark preceding Proposition \ref{p2}).

\section*{4. Galois groups and irreducibility}

As in the foregoing section, let $n\ge 3$ be odd, $d,R\in\Q$, $d\ne 0$, $R$ not a square, $D=d^2-R$, and $f_n=f_n(Z,d,R)$.
Further, $L$ denotes the splitting field of $f_n$, $\zeta_n=\zeta$ a primitive $n$th root of unity and $K=\Q(\sqrt R(\zeta-\zeta^{-1}))$.

\begin{prop} 
\label{p4}

The group $\GAL(K/\Q)$ is isomorphic to $(\Z/n\Z)^{\times}$ except in the case $R<0, \sqrt R\in\Q(\zeta)$, where this group is isomorphic to $(\Z/n\Z)^{\times}/\{\pm 1\}$.

\end{prop} 

\MN
{\em Proof.} Suppose, first, that $\sqrt R\in\Q(\zeta)$. Then $K\subseteq \Q(\zeta)$. If $R<0$, then $\sqrt R(\zeta-\zeta^{-1})$ is real and, thereby, $K\subseteq\Q(\zeta)^+$ (recall the proof of Proposition \ref{p1}).
But $K$ contains $\Q(\zeta)^+$. Hence $K=\Q(\zeta)^+$ and $\GAL(K/\Q)\cong(\Z/n\Z)^{\times}/\{\pm 1\}$. If $R>0$, then $\Q(\zeta)^+$ is a proper subfield of $K$, and so $K=\Q(\zeta)$. This implies
$\GAL(K/\Q)\cong(\Z/n\Z)^{\times}$.

Now suppose that $\sqrt R\not\in\Q(\zeta)$. Then $K$ is a subfield of $\Q(\sqrt R,\zeta)$. For an integer $k$, $(k,n)=1$, let $\sigma_k\in\GAL(\Q(\sqrt R,\zeta)/\Q)$ be defined by $\sigma_k(\zeta)=\zeta^k$,
$\sigma_k(\sqrt R)=\sqrt R$.
The subgroup $G=\{\sigma_k; (k,n)=1\}$ of $\GAL(\Q(\sqrt R,\zeta)/\Q)$ is isomorphic to $\GAL(\Q(\zeta)/\Q)$ and, thus, to $(\Z/n\Z)^{\times}$. Since $\Q(\zeta)=\Q(\zeta-\zeta^{-1})$, the numbers
$\zeta^k-\zeta^{-k}$, $0\le k\le n-1$, $(k,n)=1$, are distinct. Therefore, the map $G\to \GAL(K/\Q):\sigma_k\mapsto\sigma_k|_K$ is a monomorphism. It is surjective, since for $\tau\in\GAL(\Q(\sqrt R,\zeta)/\Q)$
defined by $\tau(\sqrt R)=-\sqrt R$, $\tau(\zeta)=\zeta$, we have $\tau|_K=\sigma_{-1}|_K$.
\STOP

\MN
As above, let $u=z^{(n-1)/2}(y+y')$, where $y$, $y'$ are zeros of $h,h'$, respectively, and $z=yy'$. We put
\BD
\label{4.2}
   x=z^{(n-1)/2}y
\ED
and $K_0=\Q(\sqrt R)$. Then $x$ is a zero of $P=Z^n-D^{(n-1)/2}(d+\sqrt R)\in K_0[Z]$, $D=d^2-R$. We show the following.

\begin{lemma} 
\label{l1}

The polynomial $f_n\in\Q[Z]$ is irreducible if, and only if, $P\in K_0[Z]$ is irreducible.

\end{lemma} 

\MN
{\em Proof.} We argue with the degrees of fields over subfields.
If $f_n$ is irreducible, then $[\Q(u):\Q]=n$. Since $n$ is odd and $[K_0:\Q]=2$, we obtain $[K_0(u):K_0]=n$. By
\BE
\label{4.4}
 u=x+D/x,
\EE
$K_0(u)\subseteq K_0(x)$. But $[K_0(x):K_0]\le n$, and so
$[K_0(x):K_0]=n$ and $P$ is irreducible. Conversely, if $[K_0(x):K_0]=n$, we obtain $[K_0(u):K_0]\in\{n,n/2\}$ because $x$ satisfies the quadratic equation (\ref{4.4}) over $K_0(u)$. However, $n$ is odd.
Accordingly, $[K_0(u):K_0]=n$ and, thus, $[\Q(u):\Q]=n$.
\STOP

\MN
The {\em proof of Theorem \ref{t1}} is contained in the following two propositions.

\begin{prop} 
\label{p5}
Let $f_n\in\Q[Z]$ be irreducible and $L$ the splitting field of $f_n$ over $\Q$. Let $x$ be defined as above, in particular, $x$ is a zero of $P$. Let $\zeta$ denote a primitive $n$th root of unity.
Suppose that $K_0(x)\cap K_0(\zeta)=K_0$.
Then the assertion of Theorem \ref{t1} holds for $\GAL(L/\Q)$.

\end{prop} 

\MN
{\em Proof.} Obviously, $K_0(x,\zeta)$ is the splitting field of $P$ over $K_0$, and, since $u\in K_0(x)$, $L=\Q(u)\cdot K\subseteq K_0(x,\zeta)$.
We have
\BD
 \GAL(K_0(x,\zeta)/K_0(x))\cong\GAL(K_0(\zeta)/K_0(x)\cap K_0(\zeta)),
\ED
see \cite[p. 196]{La}. By assumption, $K_0(x)\cap K_0(\zeta)=K_0$, and so $\GAL(K_0(x,\zeta)/K_0(x))\cong\GAL(K_0(\zeta)/K_0)$. In particular,
$[K_0(x,\zeta):K_0(x)]=[K_0(\zeta):K_0]$. By Lemma \ref{l1}, $[K_0(x):K_0]=n$. From the identities
\BD
 [K_0(x,\zeta):K_0]=[K_0(x,\zeta):K_0(x)][K_0(x):K_0]=[K_0(\zeta):K_0]\cdot n
\ED
and
\BD
  [K_0(x,\zeta):K_0]=[K_0(x,\zeta):K_0(\zeta)][K_0(\zeta):K_0]
\ED
we conclude $[K_0(x,\zeta):K_0(\zeta)]=n$. Accordingly, $P$ is irreducible over $K_0(\zeta)$. This means that $\GAL(K_0(x,\zeta)/K_0(\zeta))$
acts on the zeros of $P$ transitively. We denote these zeros by
\BD
  x_k=x\zeta^k,\: k\in\Z/n\Z
\ED
(see the remark at the beginning of Section 3). Hence there is a
$\sigma\in\GAL(K_0(x,\zeta)/K_0(\zeta))$ with $\sigma(x)=x_1$. This implies $\sigma(x_k)=\sigma(x\zeta^k)=x_1\zeta^k=x_{k+1}$ for all $k\in\Z/n\Z$. Thus $\sigma$ has the order $n$.
Now (\ref{4.4}) shows that $\sigma(u_k)=u_{k+1}$ for all $k\in\Z/nZ$. Since $\sigma$ fixes $\zeta$ and $\sqrt R$, we have $\sigma|_L\in\GAL(L/K)$. But $\sigma|_L$ also has the order $n$,
and $[L:K]=[K(u):K]\le n$. So we see that $\sigma$ generates $\GAL(L/K)$, a cyclic group of order $n$.
Because $K$ is an abelian extension of $\Q$, $\GAL(L/K)$ is a normal subgroup of $\GAL(L/\Q)$.

On the other hand, $\GAL(L/\Q(u))\cong\GAL(K/K\cap\Q(u))$. But $K\subseteq K_0(\zeta)$ and $\Q(u)\subseteq K_0(x)$, hence $K\cap\Q(u)\subseteq K_0(\zeta)\cap K_0(x)=K_0$. However, $\Q(u)\cap K_0=\Q$.
Altogether, $\GAL(L/\Q(u))\cong\GAL(K/\Q)$. By Proposition \ref{p4}, $\GAL(L/\Q(u))$ is isomorphic to $(\Z/n\Z)^{\times}$ except in the case $R<0, \sqrt R\in\Q(\zeta)$,
where this group is isomorphic to $(\Z/n\Z)^{\times}/\{\pm 1\}$.
By Galois theory, we have $\GAL(L/K)\cap\GAL(L/\Q(u))=\GAL (L/K\cdot \Q(u))=\{id\}$. Hence $\GAL(L/\Q)$ contains the semidirect product $\GAL(L/K)\rtimes\GAL(L/\Q(u))$, and an inspection of the degrees of these groups
shows that $\GAL(L/\Q)$ equals this semidirect product.
\STOP

\MN
Recall that $R$ has the form $R=s^2R'$ with $s\in\Q$ and $R'\in \Z$ square-free.

\begin{prop} 
\label{p6}
Let $f_n\in\Q[Z]$ be irreducible and $x$ be defined as above, in particular, $x$ is a zero of $P$. Let $\zeta$ denote a primitive $n$th root of unity.
If $K_0(x)\cap K_0(\zeta)\ne K_0$, then $3\DIV n$ and $R'=-3$.

\end{prop} 

\MN
{\em Proof.} Suppose $K_0(x)\cap K_0(\zeta)\ne K_0$. Then $K_0(x)\cap K_0(\zeta)$ contains a minimal subfield of $K_0(x)$ different from $K_0$. By \cite[Thm. 2.3]{OrVe}, this field has the form
$K_0(w)$, for some number $w=x^{n/p}\zeta_p^j$, where $p\DIV n$, $\zeta_p$ is a primitive $p$th root of unity, and $j\in\Z$. Since $K_0(w)$ is contained in $K_0(\zeta)$, it is an abelian extension of $K_0$.
By Lemma \ref{l1}, $P$ is irreducible over $K_0$ and, thus, $Z^p-D^{(n-1)/2}(d+\sqrt R)$ is irreducible, see \cite[Th. 1.2]{MoVe}. Since $w$ is a zero of this polynomial, we have
$[K_0(w):K_0]=p$. Accordingly,
the number $w\zeta_p$ is conjugate to $w$. Since $K_0(w)$ is abelian over $K_0$, this conjugate is contained in $K_0(w)$ and, thereby, $\zeta_p\in K_0(w)$. We have $[\Q(\zeta_p):\Q]=p-1$ and, thus,
$[K_0(\zeta_p]:K_0]=(p-1)/2$ if $\sqrt R\in \Q(\zeta_p)$ and $[K_0(\zeta_p]:K_0]=p-1$, otherwise. Therefore, $\zeta_p\in K_0(w)$ implies that $(p-1)/2$ divides $p$. This can be the case only if $p=3$ and $\sqrt R\in\Q(\zeta_3)$,
i.e., $R'=-3$.
\STOP

\MN{\em Remarks.}
1. Propositions \ref{p5} and \ref{p6} immediately yield Theorem \ref{t1}.

2. Suppose that $3\DIV n$ but $9\NDIV n$. Then the proof of Proposition \ref{p6} shows that $K_0(x)\cap K_0(\zeta)\ne K_0$ only if $R'=-3$ and $n$ is divisible by a prime $p\equiv 1$ mod 3.
Indeed, we have  $K_0(x)\cap K_0(\zeta)\ne K_0$ if $n=21$, $d=7/2$, and $R=21^2(-3)/4$, for instance.

3. Theorem \ref{t1} combined with the last remark shows that the assertion of Theorem \ref{t1} holds if $n=p$ is a prime $\ge 3$.

4. In the case $n=3^r$, $r\ge 2$, and $R'=-3$, one can show that $K_0(x)\cap K_0(\zeta)\ne K_0$ except if $P=Z^n-a^3\zeta_3$ or $P=Z^n-a^3\zeta_3^2$ for some $a\in K_0$.

\MN
{\em Proof of Theorem \ref{t2}.} By Lemma \ref{l1}, $f_n$ is irreducible, if, and only if, $P$ is irreducible in $K_0[Z]$. Let $x$ be a zero of $P$. Then it is easy to check
that $x_1=x^2/D$ is a zero of $P_1=Z^n-(d+\sqrt R)/(d-\sqrt R)$. Conversely, if $x_1$ is a zero of $P_1$, then $x=x_1^{(n+1)/2}(d-\sqrt R)$ is a zero of $P$.
Therefore, $P$ is irreducible if, and only if, $P_1$ is irreducible. Now the theorem follows from the irreducibility criterion \cite[Th. 1.2]{MoVe} applied to $P_1$.
\STOP

\MN
Theorem \ref{t2} yields sufficient conditions for the irreducibility of $f_n$ that can be checked easily.
To this end we use the ring of integers $\OO$ of $K_0=\Q(\sqrt R)$ and the valuation $v_q:\Q\sminus\{0\}\to\Z$ for a prime number $q$.
We also use the valuation $v_Q:K_0\sminus\{0\}\to\Z$ for a prime ideal $Q$ in $\OO$.

\begin{prop} 
\label{p7}
Let $d$, $R$ be integers, $R$ not a square, with $(d,R)=1$ and $D=d^2-R$.  Suppose that for each prime divisor $p$ of $n$ there is a prime $q\ge 3$ such that
$v_q(D)$ is odd and $\not\equiv 0$ mod $p$. Then $f_n=f_n(Z,d,R)$ is irreducible.

\end{prop} 

\MN
{\em Proof.} Let $p$ be a prime divisor of $n$ and $q$ as in the proposition.
First we show that $q$ is split in $\OO$. Indeed, if $q$ is ramified, it divides $R$ and, thus, the number $d$, which is excluded by $(d,R)=1$.
If $q$ is inert, then $v_q(d+\sqrt R)=v_q(d-\sqrt R)$ and $v_q(d^2-R)=2v_q(d+\sqrt R)$ is even.
So we have $\SP{q}=QQ'$, where $Q$, $Q'$ are distinct prime ideals in $\OO$.
Since $Q$ divides $d^2-R$, we may assume, without loss of generality, that $v_Q(d+\sqrt R)>0$. Then $v_Q(d-\sqrt R)=0$, since, otherwise, $Q$ divides the ideal
$\SP{d+\sqrt R,d-\sqrt R}$, which, in turn, divides $\SP{2}$, a contradiction. Accordingly, $v_Q((d+\sqrt R)/(d-\sqrt R))=v_Q(d+\sqrt R)=v_Q(D)=v_q(D)\not\equiv 0$ mod $p$.
In particular, $(d+\sqrt R)/(d-\sqrt R)$ is not a $p$th power in $K_0$.
By Theorem \ref{t2}, $f_n$ is irreducible.
\STOP

\MN
Of course, Proposition \ref{p7} yields numerous examples of irreducible De Moivre polynomials $f_n$.
The next theorem is contained in the paper \cite{SpWi} in the special case $n=5$.
As above, let $n\ge 3$ be odd, $d,R\in\Q$, $d\ne 0$, $R$ not a square, $D=d^2-R$, and $f_n=f_n(Z,d,R)$.

\begin{theorem} 
\label{t4}
If $n=p$ is a prime $\ge 3$ and $f_p$ has no rational zero, then $f_p$ is irreducible in $\Q[Z]$.

\end{theorem} 

\MN
{\em Proof.}
Suppose that $f_p$ is reducible. According to Lemma \ref{l1}, $P=Z^p-D^{(p-1)/2}(d+\sqrt R)$ is reducible over $K_0=\Q(\sqrt R)$. By \cite[Th. 1.2]{MoVe}, $P$ has a zero of the form $a+b\sqrt R$, $a,b\in\Q$.
As above, let $x=z^{(p-1)/2}y$. Then $x$ is a zero of $P$, i.e. $x=(a+b\sqrt R)\zeta^k$ for a primitive $p$th root of unity $\zeta$ and $k\in\Z$.
We may assume that $k=0$. Otherwise, we replace $y$ by $y\zeta^{-k}$ and $y'$ by
$y'\zeta^k$. From (\ref{4.4}) we conclude that $u=x+D/x$ is a zero of $f_p$. Now
$(a+b\sqrt R)^p=D^{(p-1)/2}(d+\sqrt R)$ implies $(a-b\sqrt R)^p=D^{(p-1)/2}(d-\sqrt R)$. These identities yield
\BD
 (a^2-b^2R)^p=(a+\sqrt R)^p(a-\sqrt R)^p=D^{p-1}(d+\sqrt R)(d-\sqrt R)=D^p.
\ED
However, $D^p$ has only one $p$th root in $\Q$, namely, the number $D$. This implies $a^2-b^2R=D$, and, thus, $D/x=a-b\sqrt R$. Hence $f_p$ has the zero $a+b\sqrt R+a-b\sqrt R=2a\in\Q$.
\STOP

\MN
{\em Remark.} If $n=p$ is a prime, Theorem \ref{t4}, Theorem \ref{t1}, Corollary \ref{c1}, and the third remark following Proposition \ref{p6} show that the structure of $\GAL(L/\Q)$ is completely known in all cases.

\MN
We use Theorem \ref{t4} in order to obtain a result of Filaseta et al. in a simple way, see \cite{Fi}. These authors consider the polynomial $f_p$ of prime degree $p\ge 5$ in the case
$d=-1/(2m^{(p-1)/2})$ and $R=(1-4m^p)/(4m^{p-1})$ for a natural number $m\ge 2$. Of course, $d\ne 0$ and $R$ is negative, in particular, not a square. By (\ref{1.12}), the polynomial $f_p$ takes the simple form
\BD
 f_p=\sqrt m^{\:p}F_p(Z/\sqrt m)+1=1+\sum_{k=0}^{(p-1)/2}(-1)^k\frac p{p-k}\binom{p-k}{k}m^kZ^{p-2k}.
\ED
A considerable part of the aforesaid paper is devoted to the proof of the irreducibility of $f_p$. We can see this as follows.
If $f_p$ is reducible, Theorem \ref{t4} says that it has a zero $u$ in $\Q$. However, $f_p$ is a monic polynomial with integer coefficients, and $f_p(0)=1$. Therefore, $u$ can take only the values $\pm 1$.
We have
\BD
 f_p(1)=2+\sum_{k=1}^{(p-1)/2}(-1)^k\frac p{p-k}\binom{p-k}{k}m^k.
\ED
The sum on the right-hand side of this formula is divisible by $p$, and so $f_p(1)\ne 0$. In the remaining case, we have
\BD
  f_p(-1)=\sum_{k=1}^{(p-1)/2}\frac p{p-k}\binom{p-k}{k}m^k(-1)^{p-k}.
\ED
Suppose that $f_p(-1)=0$.
If we divide by $m$, we obtain
\BD
  p+\sum_{k=2}^{(p-1)/2}\frac p{p-k}\binom{p-k}{k}m^{k-1}(-1)^{p-k}=0
\ED
(recall $p\ge 5)$. But the sum on the left side is divisible by $m$, and so $m$ divides $p$. Since $m\ge 2$, we have $m=p$.
Then the aforesaid sum is divisible by $p^2$, and so $p^2$ divides $p$. Hence $f_p(-1)=0$ is impossible.

\section*{5. Remarks on the case: $n$ even }

In this section we collect some facts about the case of an even degree $n\ge 2$ without proofs.
Let $d,R\in\Q$, $d\ne 0$, $R$ not a square, and $h,h'$ be defined as in (\ref{1.9}). Let $y,y'$ be zeros of $h,h'$, respectively. In general, $y$, and $y'$ are irrationals of degree $2n$.
A reduction of $y,y'$ to irrationals of degree $<2n$ is only possible in special cases. Consider, for instance, the case $n=4$, $d=2$, $R=7$. Here $y$ is a zero of the irreducible polynomial $g=Z^8-4Z^4-3$.
Let $M$ be the splitting field of $g$ over $\Q$. Then $\GAL(M/\Q)$ contains a cycle of length 8, which implies that $M$ cannot be a subfield of a composite $L_1\cdots L_m$, where $L_j$ is the splitting field of a
polynomial of degree $\le 7$, $j=1,\LD,m$. Therefore, a reduction of this kind is impossible.

However, if $R=d^2-v^n$ for some $v$ in $\Q$, then it is possible to reduce $y,y'$ by means of a De Moivre Polynomial. It is defined by (\ref{1.16}) and has the zero $u=y+y'$, where $y$, $y'$ have been chosen such
that $yy'=v$.
We write $f_n=f_n(Z,d,v)$ for this polynomial.

Instead of (\ref{2.2}), we have
the reduction formula
\BD
  \{y,y'\}=\left\{ u/2 \pm A(u)\sqrt{R}/2\right\},
\ED
where $A$ is the polynomial $(\sqrt v^{\:n+1}F_{n-1}(Z/\sqrt v)-dZ)/R$. Let $u_k=y\zeta^k+y'\zeta^{-k}$, $k=0,\LD,n-1$, be the zeros of $f_n$. Then the reduction formula gives
\BD
 u_k=(u/2)(\zeta^k+\zeta^{-k})+(A(u)/2)\sqrt R(\zeta^k-\zeta^{-k}),
\ED
$k=0,\LD,n-1$. Thereby, the splitting field $L$ of $f_n$ becomes
\BD
  L=\Q(u,\zeta+\zeta^{-1},\sqrt R(\zeta-\zeta^{-1})).
\ED
Note that $\zeta+\zeta^{-1}$ cannot be omitted, since $(\sqrt R(\zeta-\zeta^{-1}))^2\in\Q(\zeta^2+\zeta^{-2})$, which is a proper subfield of $\Q( \zeta+\zeta^{-1})$.

Suppose we want to obtain an analogue of Theorem \ref{t1} in the present case.
The role of $x$ is played by the zero $y$ of $h$. But Lemma \ref{l1}, with $y$ instead of $x$, is no longer true. Indeed, it may happen that $f_n$ is irreducible over $\Q$
but $h$ is reducible over $K_0=\Q(\sqrt R)$. This happens, for instance, if $n=4$, $d=3$, $v=1$. Conversely, if $h$ is irreducible over $K_0$, then $f_n$ is irreducible over $\Q$.
Accordingly, we must assume that $h$ is irreducible over $K_0$. An irreducibility criterion for $h$ has to consider the case $p=2$ additionally. Here it is not sufficient, in general,
that $d+\sqrt R$ is not a square in $K_0$. In fact, if $4\DIV n$, we also need that $-4(d+\sqrt R)$ is not a fourth power in $K_0$, see \cite[Th. 1.2]{MoVe}.

The condition $K_0(y)\cap K_0(\zeta)=K_0$, $\zeta=\zeta_n$,
is needed for an analogue of Proposition \ref{p5}. We briefly highlight the case $n=2^r$, $r\ge 3$.
Suppose that $R>0$ and $d+\sqrt R>0$. Then $i\not\in K_0(y)$, since $K_0(y)$ is conjugate to a real field. By
 \cite[Thm. 2.1]{OrVe}, $K_0(y)$ has a smallest subfield $\ne K_0$ over $K_0$, namely, $K_0(y^{n/2})$. Hence, if $K_0(y)\cap K_0(\zeta)\ne K_0$,
the field $K_0(y^{n/2})$ must coincide with a quadratic extension of $K_0$ in $K_0(\zeta)$. Suppose that $R'\ne 2$. Since $K_0$ does not contain $\sqrt 2$, $i$, and $i\sqrt 2$,
the only quadratic extensions of $K_0$ in $K_0(\zeta)$ are $K_0(\sqrt 2)$, $K_0(i)$, and $K_0(i\sqrt 2)$. Of course, $K_0(y^{n/2})$ can coincide only with $K_0(\sqrt 2)$. This, however,
does not happen unless $d$ has the form $d=4w^2\pm v^{n/2}$ for some $w\in\Q$.
Under these conditions, we obtain $\GAL(L/\Q)\cong(\Z/nZ)\rtimes(\Z/nZ)^{\times}$. An example of this kind is given by $n=8, d=4,v=1$.

The conditions $R>0$ and $d+\sqrt R>0$ have been chosen for the sake of simplicity. They are certainly not best possible.
The above discussion shows that we are faced not only with the restriction $R=d^2-v^n$, $v\in\Q$, but also with new complications in the case of an even degree $n$.


\MN
Kurt Girstmair\\
Institut f\"ur Mathematik \\
Universit\"at Innsbruck   \\
Technikerstr. 13/7        \\
A-6020 Innsbruck, Austria \\
Kurt.Girstmair@uibk.ac.at

\end{document}